\documentclass[a4paper, 10pt]{article}

\usepackage{cite}

\usepackage[top=1.5in,bottom=1.5in,right=1.5in,left=1.5in]{geometry}
\usepackage{amsfonts}
\usepackage{amsmath}
\usepackage{amssymb}
\usepackage{amsthm}
\usepackage{fancyref}
\usepackage{mathrsfs}
\usepackage{enumitem}
\usepackage{color}
\usepackage{romannum}
\usepackage{mathtools}
\usepackage{csquotes}
\usepackage[UKenglish]{babel}
\usepackage[UKenglish]{isodate}
\usepackage[capitalise]{cleveref}
\usepackage[normalem]{ulem}
  \crefformat{equation}{\textup{#2(#1)#3}}

\theoremstyle{plain}
\newtheorem{theorem}{Theorem}[section]
\newtheorem{lemma}[theorem]{Lemma}
\newtheorem{proposition}[theorem]{Proposition}
\newtheorem{corollary}[theorem]{Corollary}
\theoremstyle{remark}
\newtheorem{remark}[theorem]{Remark}
\theoremstyle{definition}
\newtheorem{definition}[theorem]{Definition}
\newtheorem{example}[theorem]{Example}
\newtheorem{assumption}[theorem]{Assumption}

\crefname{assumption}{Assumption}{Assumptions}
\Crefname{assumption}{Assumption}{Assumptions}

\newcounter{counter_a}
\newenvironment{myenum}{\begin{list}{\textrm{\textup{(\roman{counter_a})}}}%
{\usecounter{counter_a}
\setlength{\itemsep}{0.5ex}\setlength{\topsep}{0.7ex}
\setlength{\leftmargin}{5ex}\setlength{\labelwidth}{5ex}}}{\end{list}}

\newcounter{counter_b}
\newenvironment{myenuma}{\begin{list}{\textrm{\textup{(\alph{counter_b})}}}%
{\usecounter{counter_b}
\setlength{\itemsep}{0.5ex}\setlength{\topsep}{0.7ex}
\setlength{\leftmargin}{5ex}\setlength{\labelwidth}{5ex}}}{\end{list}}

\DeclareMathOperator\real{Re}
\DeclareMathOperator\imag{Im}
\renewcommand\Re{\real}
\renewcommand\Im{\imag}

\DeclareMathOperator*\esssup{ess\,sup}
\newcommand\rd{\mathrm{d}}
\newcommand\CC{\mathbb{C}}
\newcommand\NN{\mathbb{N}}
\newcommand\RR{\mathbb{R}}

\newcommand\cD{\mathcal{D}}

\newcommand\normcdot{\lVert\,\cdot\,\rVert}
\newcommand\normcdotsub[1]{\lVert\,\cdot\,\rVert_{#1}}

\newcommand{\bm}{b_{\textup{\textsf{m}}}}
\newcommand{\BigO}{\mathrm{O}}

\newcommand{\Linfloc}{L^\infty_{\textup{\textsf{loc}}}}
\definecolor{grey}{rgb}{0.4,0.4,0.4}
\definecolor{darkgreen}{rgb}{0,0.6,0}
\definecolor{MLcol}{rgb}{0.5,0,0.8}

\numberwithin{equation}{section}

\begin{document}

\pagenumbering{arabic}

\title{Continuous Fragmentation Equations in \\[0.5ex] Weighted $L^1$ Spaces \\[2ex] 
{\normalsize \textit{Dedicated to Rainer Picard on the occasion of his 80th birthday}}}
\date{}

\author{Lyndsay Kerr, Wilson Lamb and Matthias Langer}

\maketitle

\begin{abstract}
\noindent
We investigate an integro-differential equation that models the evolution of fragmenting clusters. 
We assume cluster size to be a continuous variable and allow for situations in which mass is not necessarily 
conserved during each fragmentation event.  
We formulate the initial-value problem as an abstract Cauchy problem (ACP) in an appropriate weighted $L^1$ space, 
and apply perturbation results to prove that a unique, physically relevant classical solution of the ACP 
is given by a strongly continuous semigroup for a wide class of initial conditions. 
Moreover, we show that it is often possible to identify a weighted $L^1$ space in which this semigroup is analytic, 
leading to the existence of a unique, physically relevant classical solution for all initial conditions belonging to that space. 
For some specific fragmentation coefficients, we provide examples of weighted $L^1$ spaces where our results can be applied.
\\[1ex]
\textit{Mathematics Subject Classification \textup{(}2020\textup{)}: 47D06, 80A30, 46B42}
\\[1ex]
\end{abstract}


\section{Introduction}
\label{E11}

There are many processes that can be described in terms of clusters of particles that can merge together (i.e.\ \emph{coagulate}) 
to produce larger clusters and can break apart (i.e.\ \emph{fragment}) to produce smaller clusters. 
As highlighted in \cite[Section~2.1]{banasiak_lamb_laurencot2020vol1}, such processes arising in nature include 
the grouping of animals and fish, the formation of preplanetisimals and blood clotting. 
Coagulation and fragmentation processes also feature in a number of industrial applications, and play an important role 
in the powder production industry \cite{verdurmen_etal2004simulation, wells2018thesis}, aerosols \cite{drake1972aerosol} 
and the formation and degradation of polymers \cite{aizenman_bak1979, ziff1980kinetics, ziff_mcgrady1985kinetics}.  
When clusters are comprised of identical particles (monomers), then the evolution of clusters can be described 
in terms of an infinite system of ordinary differential equations. 
We have previously examined such discrete models in the case of pure fragmentation, both time-dependent and time-independent, 
\cite{kerr_lamb_langer2020fragpaper, kerr_lamb_langer2024nonautonomous} and in the case of coagulation--fragmentation 
with time-dependent coagulation \cite{kerr_langer2025discrete}. 

In this article, we turn our attention to continuous models where cluster mass, which we use interchangeably with cluster size, 
can take any positive value. We focus on the case where no coagulation occurs, i.e.\ the breakdown of clusters is non-reversible. 
In this case, the evolution of clusters can be described in terms of the integro-differential equation  
\begin{equation}\label{E40}
\begin{split}
	\frac{\partial}{\partial t}u(x,t) &= -a(x)u(x,t)+\int_x^{\infty} a(y)b(x,y)u(y,t) \, \rd y, \qquad x,t>0; 
	\\[1ex]
	u(x,0) &= \mathring{u}(x), \qquad x>0,
\end{split}
\end{equation}
where $u(x,t)$ is the density of clusters of size $x>0$ at time $t \ge 0$, $a(x)$ is the rate of fragmentation of clusters of size $x$, 
$b(x,y)$ gives the average number of clusters of size $x$ that are produced when a cluster of size $y$ fragments 
and $\mathring{u}(x)$ is the initial density of clusters of size $x$ at time $t=0$. 
We adopt the natural physical assumption that $a(x)$ and $b(x,y)$ are non-negative for all $x,y>0$. 
We also assume that $b(x,y)=0$ whenever $y \le x$ so that a cluster cannot fragment to produce clusters larger than itself. 
Note that, since $u$ represents a density, any physically relevant solution to this system must be non-negative. 
It is also common to impose the assumption
\begin{equation}\label{E41}
	\int_0^y b(x,y)x\,\rd x = y \quad \text{for each} \ y > 0,
\end{equation}
which ensures that mass is conserved during each fragmentation event; 
however, in most of our results we do not assume that \eqref{E41} holds.

In this article, we use an operator semigroup approach to examine \eqref{E40}. 
Operator semigroups were first used, in \cite{aizenman_bak1979}, to study a binary fragmentation version of the 
continuous coagulation--fragmentation system, and have since been deployed in numerous analyses of both 
discrete and continuous coagulation--fragmentation systems; 
see \cite{banasiak2001extension, banasiak_lamb_langer2013strong, kerr_lamb_langer2020fragpaper, kerr_lamb_langer2024nonautonomous, kerr_langer2025discrete, mcbride_smith_lamb2010strongly, mclaughlin_lamb_mcbride1997fragmentation, mclaughlin_lamb_mcbride1997CF, smith_lamb_langer_mcbride2012discrete} 
and also \cite{banasiak_lamb_laurencot2020vol1}.  The approach that we use involves formulating \eqref{E40} as a linear abstract Cauchy problem (ACP) 
in an appropriate Banach space.  We then use perturbation theorems for operator semigroups to show that, 
under mild conditions on the fragmentation coefficients, there exists a unique solution of the ACP,
expressible in terms of a substochastic semigroup (i.e.\ a positive $C_0$-semigroup of contractions, often referred to as the fragmentation semigroup) 
for a particular class of initial conditions. 

For a non-negative solution, $u$, of the continuous fragmentation equation \eqref{E40}, the total number of clusters and the total mass 
of clusters at time $t \geq 0$ are given, respectively, by 
\[
	M_0(u(\cdot,t)) \coloneq \int_0^{\infty} u(x,t) \, \rd x \qquad \text{and} \qquad M_1(u(\cdot,t)) \coloneq \int_0^{\infty} u(x,t)x \, \rd x.
\]
To control both the number of particles and the mass of the system, a natural Banach space to work in is therefore the Banach space $X_{[1]}$, 
where $X_{[1]}$ is the weighted $L^1$ space,
\[
	X_{[1]} \coloneq L^1\bigl((0,\infty), (1+x)\,\rd x\bigr),
\]
consisting of real-valued, measurable functions defined (almost everywhere) on $(0,\infty)$, and satisfying 
\[
	\|f\|_{[1]} \coloneq \int_0^{\infty} |f(x)|(1+x)\,\rd x < \infty.
\]
However, while the majority of previous investigations into \eqref{E40} have utilised $X_{[1]}$, some have been conducted in other weighted $L^1$ spaces, 
primarily higher moment spaces of the form
\[
	X_{[p]} \coloneq L^1\bigl((0,\infty), (1+x^p)\,\rd x\bigr),
\]
with associated norm
\[
	\|f\|_{[p]} \coloneq \int_0^{\infty} |f(x)|(1+x^p)\,\rd x
\]
for $p \ge 1$ \cite{banasiak_lamb2012analytic, banasiak_lamb_langer2013strong}. 
Motivated by our work in \cite{kerr_lamb_langer2020fragpaper, kerr_lamb_langer2024nonautonomous, kerr_langer2025discrete}, 
here we choose to analyse \eqref{E40} in the following, more general, weighted $L^1$ space
\begin{equation}\label{E317}
	X_{\omega} \coloneq L^1\bigl((0,\infty), \omega(x)\,\rd x\bigr),
\end{equation}
where $\omega$ is a non-negative measurable function on $[0,\infty)$.  It follows that the space $X_{\omega}$ is a Banach space with norm
\begin{equation}\label{E301}
	\|f\|_{\omega} \coloneq \int_0^{\infty} |f(x)|\omega(x)\,\rd x < \infty,
\end{equation}
and consists of real-valued, measurable functions defined (almost everywhere) on $(0,\infty)$, and satisfying $\|f\|_{\omega} < \infty$.
As in \cite{kerr_lamb_langer2020fragpaper}, working with a more flexible weight allows us to prove results regarding the analyticity 
of the fragmentation semigroup.  This, in turn, allows us to yield results regarding the existence and uniqueness of solutions 
to the fragmentation ACP under weaker conditions on the fragmentation coefficients
than can be obtained if we restrict ourselves to working in $X_{[p]}$ for $p \ge 1$. 
In addition, we do not impose the usual mass-conservation assumption \eqref{E41} but instead impose a weaker assumption that allows, 
for example, mass to be lost during a fragmentation event. 
In our main result we prove that, under very mild assumptions on the fragmentation coefficients, it is always possible 
to find a weight $\omega$ such that the fragmentation semigroup in analytic on $X_{\omega}$. 
This result leads to the existence and uniqueness of solutions to the fragmentation ACP for all $\mathring{u} \in X_{\omega}$.

The article is structured as follows.  
In \cref{E10} we provide some definitions and preliminary results that we will need in our treatment of \eqref{E40}.  
In particular, we provide two perturbation theorems that we will apply to an ACP formulation of \eqref{E40}. 
In \cref{E12} we formulate \eqref{E40} as an ACP in $X_{\omega}$.  
Under very mild assumptions on the fragmentation coefficients we use the first perturbation theorem from \cref{E10} to show that, 
for a particular class of initial conditions, a unique solution of the fragmentation ACP is provided by a $C_0$-semigroup.
In \cref{E310} we apply the second perturbation result, under slightly stronger assumptions on the fragmentation coefficients, 
to show that this semigroup is analytic and provides a unique solution to the fragmentation ACP for all initial conditions in $X_{\omega}$. 
We then demonstrate in \cref{E13} that for a wide class of fragmentation coefficients it is possible to find a weight 
such that the results in \cref{E310} hold. Finally, in \cref{E14} we consider two specific examples of fragmentation coefficients, 
and provide particular weights which can be used in each of these cases.

\section{Preliminaries}\label{E10}

In this section, we supply preliminary material that we will require in our investigation into \eqref{E40}.  
In the first subsection, we begin by presenting some concepts and a result relating to \textit{AL}-spaces and positive semigroups. 
In the following subsection, we then provide two important perturbation theorems that we will later apply to the fragmentation ACP 
to obtain results regarding the existence and uniqueness of solutions.

\subsection{\textit{AL}-Spaces and Positive Semigroups}

Let us recall some definitions from the theory of Banach lattices and positive $C_0$-semigroups,
which will be crucial for our investigations.
For a detailed theory of general $C_0$-semigroups and positive semigroups, in particular,
we refer the reader to \cite{banasiak_arlotti2006perturbations, engel_nagel2000, batkai_kramer_rhandi2017positive}.

Let $X$ be a real \emph{ordered vector space}. 
We denote by $X_+$ the \emph{positive cone}, i.e.\ the set of all non-negative elements in $X$.
Likewise, for a subspace $D$ of $X$, the set of non-negative elements in $D$ is denoted by $D_+$.

Now suppose that $X$ is a \emph{vector lattice} (or \emph{Riesz space}), i.e.\ it is an ordered vector space 
and $\sup\{f,g\}$ exists in $X$ for all $f,g\in X$.  
Let $f \in X$.
Then $f_{\pm} \coloneqq \sup\{\pm f,0\}\in X_+$ and $|f| \coloneqq \sup\{f,-f\}\in X_+$ are well defined 
and satisfy $f=f_+-f_-$ and $|f|=f_++f_-$.
A subspace $Y$ of $X$ is called a \emph{sublattice} if it is closed under the mappings $f\mapsto f_+$, $f\mapsto f_-$.
A vector lattice is a \emph{Banach lattice} if it is a Banach space with norm $\normcdot$
and $|x|\le|y|$ implies $\|x\|\le\|y\|$ for all $x,y\in X$.
Note that in a Banach lattice the relation $\||x|\|=\|x\|$ is valid for all $x\in X$;
see \cite[(2.74)]{banasiak_arlotti2006perturbations}.

If $X$ is a Banach lattice and satisfies $\|f+g\|=\|f\|+\|g\|$
for all $f,g\in X_+$, then $X$ is called an \emph{AL-space}.
It follows from \cite[Theorems~2.64 and 2.65]{banasiak_arlotti2006perturbations} that, 
when $X$ is an \textit{AL}-space, there is a unique bounded linear functional, $\phi$, 
that extends $\normcdot$ from $X_+$ to~$X$.

Note that $X_\omega$, equipped with the norm $\normcdotsub{\omega}$,  is an \textit{AL}-space (where $f \le g \Leftrightarrow f(x) \le g(x)$ 
for almost all $x > 0$); see \cite[p.~100]{banasiak_lamb_laurencot2020vol1}.
The unique, bounded linear functional, $\phi_\omega$, that extends $\normcdot$ from $(X_\omega)_+$ to $X_\omega$ is given by
\[
	\phi_\omega(f) \coloneq \int_0^\infty f(x)\omega(x)\,\rd x, \qquad f \in X_\omega.
\]

Let us recall also some notions about operators and semigroups on ordered Banach spaces.
A linear operator $A$ on an ordered Banach space $X$ is called \emph{positive} if $Af\ge0$ for all $f\in X_+$.
A $C_0$-semigroup $(S(t))_{t\ge0}$ on an ordered Banach space is called a \emph{positive semigroup}
if $S(t)\ge0$ for all $t\ge0$.  A positive semigroup is called \emph{substochastic} if $\|S(t)\|\le1$ for all $t\ge0$
and \emph{stochastic} if $\|S(t)f\|=\|f\|$ for all $f\in X_+$ and $t\ge0$.

For later use we prove the following simple lemma.

\begin{lemma}\label{E205}
Let $U$ be a positive linear operator on a vector lattice, $X$, with lattice norm $\normcdot$.
Then $\|Uf\|\le\|U|f|\|$ for all $f \in X$.
\end{lemma}

\begin{proof}
Let $f \in X$.  Then $\pm f \le |f|$ and hence $\pm Uf \le U|f|$ by the positivity of $U$,
which, in turn, implies that $|Uf|\le U|f|$. 
Since $\normcdot$ is a lattice norm, it follows that
\[
	\|Uf\| = \big\||Uf|\big\| \le \big\|U|f|\big\|.
	\qedhere
\]
\end{proof}

When we deal with analytic semigroups, we also need the complexification of a real Banach lattice.
Let us recall the following notions (see, e.g.\ \cite[Section~2.2.5]{banasiak_arlotti2006perturbations} 
or \cite[Section~3.2.5]{banasiak_lamb_laurencot2020vol1}).
The \emph{complexification} $X_C$ of a real vector lattice $X$ is the set of pairs $(x,y) \in X \times X$, 
where we write $(x,y) \eqcolon x + iy$. 
We refer to $x$, $y$ and $x-iy$ respectively as the \emph{real part}, \emph{imaginary part} and \emph{complex adjoint} of $x+iy$.
Vector operations are extended in an obvious way, and the partial order in $X_C$ is defined by
\[
	x_1+iy_1 \le x_2 + iy_2  \quad\Leftrightarrow\quad  x_1 \le x_2 \;\; \text{and} \;\; y_1 = y_2.
\]

Note that $x+iy$ is a non-negative element in $X_C$ if and only if $x \in X_+$ and $y = 0$.
Motivated by the scalar case, we define a modulus by
\[
	|x+iy| \coloneq \sup_{\theta\in[0,2\pi]}(\cos\theta\cdot x + \sin\theta\cdot y), \qquad x,y\in X,
\]
which exists (see \cite[p.~62]{banasiak_arlotti2006perturbations}), and a norm by 
\[
	\|z\|_C \coloneq \| |z| \|, \qquad z\in X_C.
\]
Whenever we deal with analytic semigroups we consider them on the complexification~$X_C$.
Note that the complexification of $X_\omega$ is the weighted $L^1$ space
of complex-valued functions.

\subsection{Two Perturbation Theorems}

We now present two perturbation theorems, which we use later to prove the existence of semigroups associated with the 
fragmentation equation. We apply these two theorems in two different situations depending on the assumptions on the fragmentation kernel. 
First, we recall a perturbation result \cite[Proposition~2.4]{kerr_lamb_langer2020fragpaper},
which is based on \cite[Theorem~2.7]{thieme_voigt2006stochastic}.
The following theorem is essentially \cite[Proposition~2.4]{kerr_lamb_langer2020fragpaper}
but is reformulated slightly to make the assumptions clearer.

\begin{theorem}\label{E300}
Let $(X,\normcdot)$ and $(Z,\normcdotsub{Z})$ be \textit{AL}-spaces, such that
\begin{myenum}
\item 
$Z$ is a sublattice of $X$;
\item 
$Z$ is dense in $X$;
\item 
$(Z, \normcdotsub{Z})$ is continuously embedded in $(X, \normcdot)$.
\end{myenum}
Also, let $\phi$ and $\phi_Z$ be the linear extensions of $\normcdot$
from $X_+$ to $X$ and of $\normcdotsub{Z}$ from $Z_+$ to $Z$ respectively.
Let $A: \mathcal{D}(A) \to X$, $B: \mathcal{D}(B) \to X$ be linear operators in $X$
such that $\mathcal{D}(A) \subseteq \mathcal{D}(B)$.
Assume that the following conditions are satisfied:
\begin{myenuma}
\item 
$-A$ is positive;
\item 
$A$ generates a positive $C_0$-semigroup, $(T(t))_{t \geq 0}$, on $X$;
\item 
the semigroup $(T(t))_{t \ge 0}$ leaves $Z$ invariant and its restriction
to $Z$ is a \textup{(}necessarily positive\textup{)} $C_0$-semigroup
on $(Z, \normcdotsub{Z})$ \textup{(}in which case, the corresponding generator $\widetilde{A}$ is given by
\[
  \widetilde{A}f = Af \qquad \text{for all} \
  f \in \mathcal{D}(\widetilde{A}) = \bigl\{f \in \mathcal{D}(A) \cap Z: Af \in Z\bigr\});
\]
\item 
$B|_{\mathcal{D}(A)}$ is a positive linear operator;
\item 
$\phi((A+B)f) \le 0$ for all $f \in \mathcal{D}(A)_+$;
\item 
$(A+B)f \in Z$ and $\phi_{Z}((A+B)f) \le 0$ for all $f \in \mathcal{D}(\widetilde{A})_+$;
\item 
$\Vert Af \Vert \le \Vert f \Vert_Z$ for all $f \in \mathcal{D}(\widetilde{A})_+$.
\end{myenuma}
Then $G\coloneq\overline{A+B}$ generates a substochastic $C_0$-semigroup, $(S(t))_{t\ge0}$, on $X$,
and no other extension of $A+B$ generates a $C_0$-semigroup.
Moreover, the semigroup $(S(t))_{t \ge 0}$ leaves $Z$ invariant.
If $\phi((A+B)f)=0$ for all $f \in \mathcal{D}(A)_+$,
then $(S(t))_{t \ge 0}$ is stochastic.
\end{theorem}

Before we state the second perturbation theorem, let us recall the notion of relative boundedness.

\begin{definition}
\cite[Definition~III.2.1]{engel_nagel2000}
Let $X$ be a Banach space and let $A: \cD(A) \to X$ and $B:\cD(B) \to X$ be linear operators with $\cD(A) \subseteq \cD(B)$. 
Then $B$ is \emph{$A$-bounded} if there exist $\alpha,\beta\ge 0$ such that 
\begin{equation}\label{E201}
	\|Bf\| \le \alpha \|Af\|+\beta\|f\| \qquad \text{for all} \ f \in \cD(A).
\end{equation}
If $B$ is $A$-bounded, then the \emph{$A$-bound}, or \emph{relative bound}, is
\begin{equation}\label{E202}
	\alpha_0 \coloneq \inf\bigl\{\alpha\ge0: \ \text{there exists} \ \beta\ge0 \ \text{such that \eqref{E201} holds}\bigr\}.
\end{equation}
\hfill $\lozenge$
\end{definition} 

We also need the definition of a Miyadera perturbation.

\begin{definition}\label{E209}
\cite[Section~4.4]{banasiak_arlotti2006perturbations}
Let $X$ be a Banach space and let $A: \cD(A) \to X$ and $B: \cD(B) \to X$ be linear operators, where $\mathcal{D}(A) \subseteq \mathcal{D}(B) \subseteq X$.
Moreover, let $A$ be the generator of a $C_0$-semigroup, $(T(t))_{t\ge0}$, on $X$.  
Then $B$ is a \emph{Miyadera perturbation} of $A$ if $B$ is $A$-bounded and there exist numbers 
$\zeta>0$ and $\gamma\in(0,1)$ such that
\begin{equation}\label{E208}
	\int_0^{\zeta} \|BT(t)f\|\,\rd t \le \gamma \|f\| \qquad \text{for all} \ f \in \cD(A).
\end{equation}
\hfill $\lozenge$
\end{definition}

The Miyadera perturbation theorem \cite[Theorem 4.16]{banasiak_arlotti2006perturbations} states that $A+B$ generates a $C_0$-semigroup
if $A$ is the generator of a $C_0$-semigroup and $B$ is a Miyadera perturbation of $A$. 
This result is used in the proof of the following perturbation theorem.

\begin{theorem}\label{E203}
Let $A$ be the generator of a positive $C_0$-semigroup on an \textit{AL}-space, $X$, such that $-A$ is positive.  
Moreover, let $B$ be an $A$-bounded linear operator, with $A$-bound strictly less than $1$.
Then the following statements hold.
\begin{myenum}
\item
	The operator $A+B$ is the generator of a $C_0$-semigroup on $X$.
\item
	If $B$ is positive, then the semigroup generated by $A+B$ is positive.
\item
	If the semigroup generated by $A$ is analytic and $B$ is positive, 
	then the semigroup generated by $A+B$ is analytic.
\end{myenum}
\end{theorem} 

\begin{proof}
(i)
We show that $B$ is a Miyadera perturbation of $A$. 
Let $\phi$ be the unique linear extension of $\normcdot$ from $X_+$ to $X$ and let $(T(t))_{t\ge0}$ be the semigroup generated by $A$.
For $\zeta>0$ and $f\in\cD(A)_+$ we have
\begin{align}
	\int_0^{\zeta} \|AT(t)f\|\,\rd t 
	&= \int_0^{\zeta} \phi\bigl(-AT(t)f\bigr)\,\rd t
	= \phi\Bigl(-\int_0^{\zeta} AT(t)f\,\rd t\Bigr)
	\notag\\[1ex]
	&= \phi\Bigl(-\int_0^{\zeta} \frac{\rd}{\rd t}\bigl(T(t)f\bigr)\,\rd t \Bigr)
	= \phi\bigl(f-T(\zeta)f\bigr)
	\notag\\[1ex]
	&= \|f\|-\|T(\zeta)f\| 
	\le \|f\|.
	\label{E204}
\end{align}
The aim is to extend this inequality to all $f\in\cD(A)$.
To this end, define, for each $\delta>0$, the positive operator $T_{\delta}$ by 
$T_{\delta}f \coloneq \frac{1}{\delta}\int_0^{\delta} T(t)f\,\rd t$, $f\in X$.
If follows from \cite[Lemma~II.1.3\,(iii)]{engel_nagel2000} and \cite[Proposition~4.2.4\,(a)]{banasiak_lamb_laurencot2020vol1}
that $T_{\delta}$ maps $X$ into $\cD(A)$ and that $T_{\delta}f \to f$ as $\delta \to 0^+$ for all $f \in X$. 
Hence $T_{\delta}|f| \in \mathcal{D}(A)_+$ and $T_{\delta}|f| \to |f|$ as $\delta \to 0^+$ for all $f \in X$. 
Since, $-A$, $T(t)$ and $T_\delta$ are positive operators for $t\ge0$ and $\delta>0$, 
it follows from \cref{E205} that 
\[
	\|AT(t)T_{\delta}f\| 
	= \|(-A)T(t)T_{\delta}f\| 
	\le \|(-A)T(t)T_{\delta}|f|\| 
	= \|AT(t)T_{\delta}|f|\|
\]
for all $f\in X$.
Hence, for $f \in X$, we obtain from \eqref{E204} that
\begin{equation}\label{E206}
	\int_0^{\zeta} \|AT(t)T_{\delta}f\|\,\rd t 
	\le \int_0^{\zeta} \|AT(t)T_{\delta}|f|\|\,\rd t 
	\le \|T_{\delta}|f|\|.
\end{equation}
If we now let $f\in\cD(A)$, then with the help of \cite[Lemmas~II.1.3\,(ii) and (iv)]{engel_nagel2000} we obtain,
for $t\ge0$,
\begin{align*}
	AT(t)T_\delta f-AT(t)f
	&= T(t)\bigl[AT_\delta f-Af\bigr]
	= T(t)\Bigl(\frac{1}{\delta}\int_0^\delta T(s)Af\,\rd s-Af\Bigr)
	\\[1ex]
	&= T(t)\frac{1}{\delta}\int_0^\delta \bigl(T(s)-I\bigr)Af\,\rd s,
\end{align*}
which implies that
\begin{align*}
	\|AT(t)T_{\delta}f-AT(t)f\| 
	&\le \|T(t)\|\frac{1}{\delta}\int_0^{\delta} \|(T(s)-I)Af\|\,\rd s
	\\[1ex]
	&\le \|T(t)\| \sup_{s\in[0,\delta]} \|(T(s)-I)Af\|
	\to 0 \quad \text{as} \ \delta \to 0^+
\end{align*}
uniformly in $t$ on $[0,\zeta]$ for every $\zeta>0$.  
This shows that the left-hand side of \eqref{E206} converges to $\int_0^\zeta \|AT(t)f\|\,\rd t$ as $\delta\to0^+$.
Taking the limit as $\delta \to 0^+$ in \eqref{E206} we obtain
\begin{equation}\label{E207}
	\int_0^{\zeta} \|AT(t)f\|\,\rd t \le \|f\| \qquad \text{for all} \ f \in \cD(A).
\end{equation}
Let $M \ge 1$ and $\omega \ge 0$ be such that $\|T(t)\| \le Me^{\omega t}$ for all $t \ge 0$.  
By assumption there exist $\alpha\in[0,1)$ and $\beta\ge0$ such that \eqref{E201} holds.
This relation, together with \eqref{E207}, implies that
\begin{align*}
	\int_0^{\zeta} \|BT(t)f\|\,\rd t 
	&\le \alpha\int_0^{\zeta} \|AT(t)f\|\,\rd t + \beta\int_0^{\zeta} \|T(t)f\|\,\rd t
	\\[1ex]
	&\le \alpha\|f\| + \beta M\int_0^{\zeta} e^{\omega t}\|f\|\,\rd t 
	\le \bigl(\alpha+\beta M\zeta e^{\omega\zeta}\bigr)\|f\|
\end{align*}
for all $f \in \cD(A)$.
Since $\alpha<1$, we can choose $\zeta>0$ such that $\alpha+\beta M\zeta e^{\omega\zeta}<1$
and hence \eqref{E208} holds with $\gamma<1$.
Now the Miyadera perturbation theorem implies that $A+B$ generates a $C_0$-semigroup.

(ii)
Let $(T(t))_{t\ge0}$ and $(S(t))_{t\ge0}$ be the semigroups generated by $A$ and $A+B$ respectively,
and let $\omega_0$ be the growth bound of $(T(t))_{t\ge0}$. 
Choose $\lambda>\omega_0$; then $A-\lambda I$ generates the semigroup $(e^{-\lambda t}T(t))_{t\ge0}$, 
which has a negative growth bound, and $A-\lambda I+B$ generates the semigroup $(e^{-\lambda t}S(t))_{t\ge0}$. 
It follows from \cite[Lemma~4.15]{banasiak_arlotti2006perturbations} that $B$ is a Miyadera Perturbation of $A-\lambda I$.
Further, the proof of \cite[Theorem~4.16]{banasiak_arlotti2006perturbations} yields that $S(t)$ can be 
written as a series,
\begin{equation}\label{E210}
	S(t) = e^{\lambda t}\sum\limits_{j=0}^{\infty} S_j(t), \qquad t \ge 0,
\end{equation}
where $S_j(t)$, $j\in\NN_0$, $t\in[0,\infty)$, are bounded operators satisfying
\begin{alignat*}{2}
	S_0(t)f &= e^{-\lambda t}T(t)f, \qquad && f \in X; 
	\\[1ex]
	S_j(t)f &= \int_0^t S_{j-1}(t-s)Be^{-\lambda s}T(s)f \, \rd s, \qquad && f\in\cD(A), \; j=1,2,3,\ldots.
\end{alignat*}
The positivity of $(T(t))_{t\ge0}$ and $B$ imply that $S_j(t)\ge0$.
Hence the representation \eqref{E210} yields that $(S(t))_{t\ge0}$ is a positive semigroup.

(iii) 
Since $(S(t))_{t\ge0}$ is positive, the operator $A+B$ is resolvent positive; see \cite[p.~128]{banasiak_lamb_laurencot2020vol1}. 
It follows from \cite[Theorem~1.1]{arendt_rhandi1991perturbation} that $(S(t))_{t\ge0}$ is analytic.
\end{proof}

\noindent
Note that \cref{E203} could also be proved using \cite[Theorem~0.1]{voigt1989}.

\section[Fragmentation Semigroup on $X_\omega$]{Fragmentation Semigroup on {\boldmath$X_\omega$}}\label{E12}

In this section, we formulate \eqref{E40} as an ACP in $X_{\omega}$. 
We begin by introducing some natural non-negativity assumptions on the fragmentation coefficients that we require throughout our investigation.

\begin{assumption}\label{E316}
Let $a\in\Linfloc([0,\infty))$ be such that $a(x) \ge 0$ for almost all $x > 0$.
Further, let $b:(0,\infty)^2\to[0,\infty)$ be measurable such that $b(x,y)=0$ for almost all $x,y$ with $x>y$.
\hfill $\lozenge$
\end{assumption}

Let \cref{E316} hold, let $\omega:[0,\infty)\to[0,\infty)$ be measurable, and let $X_\omega$ be as in \eqref{E317}.
Motivated by \eqref{E40}, we define two operators, $A^{(\omega)}$ and $B^{(\omega)}$, in $X_\omega$ by
\[
	A^{(\omega)}f \coloneq -af, \qquad \cD(A^{(\omega)}) \coloneq \{f \in X_\omega: af \in X_\omega\},
\]
and 
\[
	(B^{(\omega)}f)(x) \coloneq \int_x^\infty a(y)b(x,y)f(y)\,\rd y, \quad  x > 0, 
	\qquad
	\cD(B^{(\omega)}) \coloneq \{f \in X_\omega: B^{(\omega)}f \in X_\omega\}.
\]
This allows us to pose \eqref{E40} as the ACP
\begin{equation}\label{E42}
	\frac{\rd u}{\rd t} =A^{(\omega)}u + B^{(\omega)}u, \quad t>0; \qquad u(0)=\mathring{u}.
\end{equation}
It follows from \cite[Propositions~I.4.10\,(i) and I.4.11]{engel_nagel2000} that the operator $A^{(\omega)}$ 
is closed, densely defined and generates the $C_0$-semigroup $(T^{(\omega)}(t))_{t\ge0}$ on $X_\omega$ given by
\[
	T^{(\omega)}(t)f = e^{-ta}f, \qquad f \in X_\omega,\, t \ge 0.
\]
This semigroup is clearly substochastic.
For later use, note that, for any $\lambda > 0$,  the resolvent operator of $A^{(\omega)}$ is given by
\begin{equation}\label{E322}
	R(\lambda,A^{(\omega)})f = \frac{1}{\lambda+a}f, \qquad f \in X_\omega,
\end{equation}
and so
\[
	\|\lambda R(\lambda,A^{(\omega)})f\|_{\omega} \le \|f\|_{\omega} \qquad 
	\text{for all} \ f \in X_\omega \ \text{and} \ \lambda > 0.
\]

We will also require the following additional assumption on the fragmentation coefficient $b$.
\begin{assumption}\label{E302}
Let $\omega:[0,\infty)\to[0,\infty)$ and assume that 
there exists $\kappa \in (0,1]$ such that
\begin{equation}\label{E303}
	n_\omega(y) \coloneq \int_0^y b(x,y)\omega(x)\,\rd x \le \kappa w(y)
	\qquad \text{for all} \ y > 0 .
\end{equation}
\hfill $\lozenge$
\end{assumption}

\medskip

\noindent
For example, suppose that
\begin{equation}\label{E304}
	\int_0^y b(x,y)x\,\rd x \le y \qquad \text{for all} \ y > 0,
\end{equation}
and also that $\gamma(x) \coloneq \omega(x)/x $ is non-decreasing on $(0,\infty)$.  
Then, for all $y > 0$, we have
\[
	\int_0^y b(x,y)\omega(x)\,\rd x = \int_0^y b(x,y)\gamma(x)x\,\rd x 
	\le \gamma(y)\int_0^y b(x,y)x\,\rd x \le \omega(y),
\]
showing that \eqref{E303} is satisfied with $\kappa = 1$. 
In particular, if \eqref{E304} holds, then \cref{E302} is automatically satisfied 
when $\omega(x) = x^p$ for any $p \ge 1$. 
We note that \eqref{E304} allows for fragmentation events in which mass is conserved or in which mass is lost. 
It follows that Assumption~\ref{E302} is weaker than the mass-conservation assumption \eqref{E41}.

In order to be able to perturb $A^{(\omega)}$ with $B^{(\omega)}$, we consider the following estimate.
It follows from \cref{E302} that, for any $f\in\cD(A^{(\omega)})_+$,
\begin{align}
	\phi_{\omega}\bigl(B^{(\omega)}f\bigr)
	&= \int_0^\infty \int_x^\infty a(y)b(x,y)f(y)\,\rd y\,\omega(x)\,\rd x
	\nonumber\\[1ex]
	&= \int_0^\infty \int_0^y b(x,y)\omega(x)\,\rd x\, a(y)f(y)\,\rd y
	\nonumber\\[1ex]
	&\le \kappa \int_0^\infty \omega(y)a(y)f(y)\,\rd y
	= -\kappa\phi_\omega\bigl(A^{(\omega)}f\bigr).
	\label{E305}
\end{align}
Consequently, for all $f \in \mathcal{D}(A^{(\omega)})$,
\begin{align}
	\|B^{(\omega)}f\|_{\omega}
	&= \int_0^\infty\bigg|\int_x^\infty a(y)b(x,y)f(y)\rd y\,\bigg|\,\omega(x)\,\rd x
	\nonumber\\[1ex]
	&\le \phi_\omega\bigl(B^{(\omega)}|f|\bigr)
	\le -\kappa\phi_\omega\bigl(A^{(\omega)}|f|\bigr)
	= \kappa\|A^{(\omega)}f\|_{\omega},
	\label{E306}
\end{align}
from which it follows that
\begin{equation}\label{E307}
	\mathcal{D}(A^{(\omega)}) \subseteq \mathcal{D}(B^{(\omega)})
	\quad\text{and}\quad
	\mathcal{D}\bigl(A^{(\omega)}+B^{(\omega)}\bigr)
	= \mathcal{D}(A^{(\omega)}) \cap \mathcal{D}(B^{(\omega)})
	= \mathcal{D}(A^{(\omega)}).
\end{equation}
To apply \cref{E300} to
the operators $A^{(\omega)}$ and $B^{(\omega)}$, we require
a suitable subspace of $X_{\omega}$.  We define such a subspace in terms of
a function $c$ that is non-decreasing on $[0,\infty)$ and satisfies
\begin{equation}\label{E318}
	a(x) \le c(x) \qquad \text{for almost every} \ x\in[0,\infty).
\end{equation}
Note that such a function can always be found when $a \in \Linfloc([0,\infty))$ as we can take
\[
	c(x) = \esssup_{y \in [0,x]}\, a(y) \qquad
	\text{for} \ x \ge 0;
\]
see \cite[Remark~5.1.38]{banasiak_lamb_laurencot2020vol1}.
Let $C^{(\omega)}$ be the multiplication operator defined by
\[
	C^{(\omega)}f = -cf, \qquad
	\cD(C^{(\omega)}) = \{f \in X_\omega: cf \in X_\omega\},
\]
and equip $\cD(C^{(\omega)})$ with the graph norm
\[
	\|f\|_{C^{(\omega)}} \coloneq \|f\|_{\omega} + \|C^{(\omega)}f\|_{\omega},
	\qquad f \in \cD(C^{(\omega)}).
\]
Note that $(\cD(C^{(\omega)}),\normcdotsub{C^{(\omega)}})$ is identical to the space 
$L^1((0,\infty),\widetilde{\omega}(x)\,\rd x) = X_{\widetilde{\omega}}$, 
where
\begin{equation}\label{E319}
	\widetilde{\omega} = (1+c)\omega.
\end{equation}

\noindent
We now apply \cref{E300} to $A^{(\omega)}$ and $B^{(\omega)}$.

\begin{theorem}\label{E308}
Let \cref{E316,E302} hold.
Then $G^{(\omega)}=\overline{A^{(\omega)}+B^{(\omega)}}$ is the generator of a
substochastic $C_0$-semigroup, $(S^{(\omega)}(t))_{t \ge 0}$, on $X_\omega$.
Moreover, if $c$ is non-decreasing and satisfies \eqref{E318} and $\widetilde\omega$ is as in \eqref{E319}, 
then $(S^{(\omega)}(t))_{t\ge0}$ leaves $X_{\widetilde{\omega}}$ invariant.
\end{theorem}

\begin{proof}
We show that the conditions of \cref{E300} are all satisfied
when $A=A^{(\omega)}$, $B=B^{(\omega)}$ and the \textit{AL}-spaces $(X,\normcdot)$
and $(Z,\normcdotsub{Z})$ are, respectively, $X_\omega$ and $X_{\widetilde{\omega}}$.

Clearly, $X_{\widetilde{\omega}}$ is a sublattice of $X_\omega$ and, furthermore, is dense and continuously embedded in $X_\omega$.
Moreover, conditions (a)--(c) are all satisfied by $A^{(\omega)}$.

It is also clear that $B^{(\omega)}$ is positive, and, for $f\in\cD(A^{(\omega)})_+$, 
\eqref{E305} leads to

\begin{equation}\label{E309}
\begin{aligned}
	\phi_\omega\bigl((A^{(\omega)}+B^{(\omega)})f\bigr)
	&= \phi_\omega(A^{(\omega)}f) + \phi_\omega(B^{(\omega)}f)
	\\[1ex]
	&\le \phi_\omega(A^{(\omega)}f) - \kappa\phi_\omega(A^{(\omega)}f) \le 0.
\end{aligned}
\end{equation}
Hence (d) and (e) hold.
Moreover, the monotonicity of the function $c$ and \cref{E302} imply that
\[
	n_{\widetilde\omega}(y) = \int_0^y b(x,y)\bigl(1+c(x)\bigr)\omega(x)\,\rd x 
	\le \bigl(1+c(y)\bigr)n_\omega(y) \le \kappa\widetilde\omega(y) \qquad \text{for all} \ y > 0.
\]
This means that \cref{E302} also holds for the weight $\widetilde\omega$.
Therefore we obtain from \eqref{E307} and \eqref{E309} that $\cD(A^{(\widetilde\omega)}) \subseteq \cD(B^{(\widetilde\omega)})$
and $\phi_{\widetilde\omega}((A^{(\widetilde\omega)}+B^{(\widetilde\omega)})f) \le 0$ 
for all $f \in \cD(A^{(\widetilde\omega)})_+$, and so (f) is also satisfied.
Finally, we use \eqref{E318} to obtain
\[
	\|A^{(\omega)}f\|_{\omega}
	\le \int_0^\infty c(x)|f(x)|\omega(x)\,\rd x
	\le \int_0^\infty |f(x)|\widetilde{\omega}(x)\,\rd x
	= \|f\|_{\widetilde{\omega}},
\]
for $f \in \mathcal{D}(\tilde{A}^{(\omega)})_+$, which shows that (g) holds.

Thus, the conditions of \cref{E300} are all satisfied
and therefore  $G^{(\omega)}=\overline{A^{(\omega)}+B^{(\omega)}}$ is the generator of
a substochastic $C_0$-semigroup, $(S^{(\omega)}(t))_{t \ge 0}$, on $X_\omega$,
which leaves $X_{\widetilde\omega}$ invariant.
\end{proof}

\begin{remark}
Note that the assumption $a\in\Linfloc([0,\infty))$ prevents the occurrence of the
phenomenon called `shattering'; see \cite[Section~2.3.1]{banasiak_lamb_laurencot2020vol1}.
\hfill $\lozenge$
\end{remark}

\medskip

\Cref{E308} allows us to deduce the solution of an ACP that is related to, but distinct from, the fragmentation ACP \eqref{E42}. 
In particular, under the conditions of \cref{E308}, $u(t)=S^{(\omega)}(t)\mathring{u}$ is the unique classical solution of the ACP
\begin{equation}\label{E43}
	\frac{\rd u}{\rd t} = G^{(\omega)}u = \overline{A^{(\omega)} + B^{(\omega)}}u, \quad t>0; \qquad u(0)=\mathring{u},
\end{equation}
for all $\mathring{u} \in \mathcal{D}(G^{(\omega)})$.  Moreover, if $\mathring{u} \in \mathcal{D}(G^{(\omega)})_+$, 
then $u(t)$ is non-negative for all $t \ge 0$.  The invariance result in \cref{E308} also allows us to obtain a solution 
of the fragmentation ACP \eqref{E42} for a certain class of initial conditions.

\begin{corollary}\label{E44}
Let \cref{E316,E302} hold.  Then $u(t)=S^{(\omega)}(t)\mathring{u}$ is the unique classical solution of \eqref{E42} 
for all $\mathring{u} \in X_{\tilde{\omega}}$. If $\mathring{u} \in (X_{\tilde{\omega}})_+$ then this solution is non-negative.
\end{corollary}

\begin{proof}
Let $\mathring{u} \in X_{\tilde{\omega}}$. From \cref{E308},  $G^{(\omega)}=\overline{A^{(\omega)}+B^{(\omega)}}$ 
is the generator of the substochastic $C_0$-semigroup, $(S^{(\omega)}(t))_{t \ge 0}$ that leaves $X_{\widetilde{\omega}}$ invariant. 
Hence, from \eqref{E318}, $S^{(\omega)}(t)\mathring{u} \in X_{\tilde{\omega}} = \mathcal{D}(C^{(\omega)}) \subseteq \mathcal{D}(A^{(\omega)})$. 
It is clear that $G^{(\omega)}$ and $A^{(\omega)}+B^{(\omega)}$ coincide on $\mathcal{D}(A^{(\omega)})$, 
and we know that $u(t)=S^{(\omega)}(t)\mathring{u}$ is the unique classical solution of \eqref{E43}. 
It follows that $u(t)=S^{(\omega)}(t)\mathring{u}$ is also the unique classical solution of \eqref{E42}. 
The non-negativity result follows from the positivity of $(S^{(\omega)}(t))_{t \ge 0}$.
\end{proof}

\begin{remark}\label{E321}
Note that $X_\omega$ is a space of type $L$, and hence there exists a measurable representation $u(x,t)$
of the classical solution that is absolutely continuous with respect to $t$ for almost all $x$
and satisfies \eqref{E40} almost everywhere; see \cite[Theorem~5.1.1]{banasiak_lamb_laurencot2020vol1}.
\hfill $\lozenge$
\end{remark}

\section{Analyticity of the Fragmentation Semigroup}\label{E310}

In this section we show that, under a slightly stronger assumption than \cref{E302}, the operator $A^{(\omega)}+B^{(\omega)}$ 
is the generator of an analytic $C_0$-semigroup.  This, in turn, allows us to deduce the existence and uniqueness of classical solutions 
to the fragmentation ACP \eqref{E42} for a wider class of initial conditions than in \cref{E44}.

\begin{assumption}\label{E211}
Let the weight $\omega: [0,\infty) \to [0,\infty)$ be such that there exist $\kappa_1 >0$, $\kappa_2 \in (0,1)$ and $\eta_0>0$ 
so that 
\begin{alignat}{2}
	\int_0^y b(x,y)\omega(x)\,\rd x &\le \kappa_1\omega(y) \qquad &&\text{for all} \ y\in(0,\eta_0],
	\label{E212}
	\\[1ex]
	\int_0^y b(x,y)\omega(x)\,\rd x &\le \kappa_2 \omega(y) \quad &&\text{for all} \ y \ge \eta_0.
	\label{E213}
\end{alignat}
\hfill $\lozenge$
\end{assumption}

\noindent
We now apply \cref{E203} to the operators $A^{(\omega)}$ and $B^{(\omega)}$.

\begin{theorem}
Let \cref{E316,E211} be satisfied.
Then the operator $G^{(\omega)}=A^{(w)}+B^{(\omega)}$ \textup{(}with domain $\cD(A^{(\omega)})$\textup{)} is the generator of 
an analytic, positive $C_0$-semigroup, $(S^{(\omega)}(t))_{t\ge0}$ on $X_\omega$.
If $\kappa_1\le1$, then $(S^{(\omega)}(t))_{t\ge0}$ is also substochastic.
\end{theorem}

\begin{proof}
We know that $A^{(\omega)}$ is the generator of a positive $C_0$-semigroup $(e^{-ta()})_{t\ge0}$ on $X_{\omega}$, 
and it is clear that $-A^{(\omega)}$ is a positive operator. 
Further, \eqref{E322} and a routine calculation shows that, for $f \in X_{\omega}$, $\lambda \in\CC\setminus\RR$ with $\Re\lambda>0$,
\[
	\big\|R(\lambda,A^{(\omega)})f\big\|_{\omega} 
	= \int_0^\infty \frac{1}{|\lambda+a(x)|}|f(x)|\omega(x)\,\rd x 
	\le \frac{1}{|\Im\lambda|}\|f\|_{\omega},
\]
and hence $(e^{-ta()})_{t\ge0}$ is an analytic semigroup by \cite[Theorem~II.4.6]{engel_nagel2000}.

We need to show that $B^{(\omega)}$ is $A^{(\omega)}$-bounded with $A^{(\omega)}$-bound strictly less than one.
Let $f \in \mathcal{D}(A^{(\omega)})$. Then
\begin{align*}
	& \|B^{(\omega)}f\|_{\omega} \le \int_0^\infty \int_x^\infty a(y)b(x,y)|f(y)|\,\rd y\,\omega(x)\,\rd x
	\\[1ex]
	&= \int_0^\infty \int_0^y b(x,y)\omega(x)\,\rd x \; a(y)|f(y)| \, \rd y
	\\[1ex]
	&= \int_0^{\eta_0} \int_0^y b(x,y)\omega(x)\,\rd x \; a(y)|f(y)|\,\rd y
	+ \int_{\eta_0}^\infty \int_0^y b(x,y)\omega(x)\,\rd x \; a(y)|f(y)|\,\rd y
	\\[1ex]
	&\le \kappa_1 \int_0^{\eta_0} \omega(y)a(y)|f(y)|\,\rd y + \kappa_2 \int_{\eta_0}^\infty \omega(y)a(y)|f(y)|\,\rd y
	\\[1ex]
	&\le \kappa_1 \esssup_{z\in[0,\eta_0]}a(z) \int_0^{\eta_0} |f(y)|\omega(y)\,\rd y 
	+ \kappa_2 \int_{\eta_0}^\infty a(y)|f(y)|\omega(y)\,\rd y
	\\[1ex]
	&\le \kappa_1 \esssup_{z\in[0,\eta_0]}a(z) \|f\|_{\omega} + \kappa_2\|A^{(\omega)}f\|_{\omega},
\end{align*}
which shows that $B^{(\omega)}$ is $A^{(\omega)}$-bounded with $A^{(\omega)}$-bound at most $\kappa_2 < 1$. 
Since $B^{(\omega)}$ is positive, it follows from \cref{E203} that $G^{(\omega)}=A^{(\omega)}+B^{(\omega)}$ 
generates a positive, analytic $C_0$-semigroup, $(S^{(\omega)}(t))_{t\ge0}$, on $X_{\omega}$. 
Moreover, if $\kappa_1\le1$, then the semigroup $(S^{(\omega)}(t))_{t\ge0}$ is substochastic by \cref{E308}.
\end{proof}

\noindent
This theorem allows us to deduce the following result regarding solutions of the fragmentation ACP.
\begin{corollary}
Let \cref{E316,E211} be satisfied.  Then $u(t)=S^{(\omega)}(t)\mathring{u}$ is the unique classical solution of \eqref{E42} 
for all $\mathring{u} \in X_{\omega}$.  If $\mathring{u} \in (X_{\omega})_+$ then this solution is non-negative. 
\end{corollary}

\begin{proof}
The existence and uniqueness follows from \cite[Theorem~1.2.4 and Corollary~4.1.5]{pazy1983semigroups}. 
The non-negativity is a consequence of the non-negativity of the semigroup $(S^{(\omega)}(t))_{t \ge 0}$.
\end{proof}

\begin{remark}\label{E216}
The fact that there exist $\eta_0>0$ and $\kappa_2<1$ such that \eqref{E213} holds is equivalent to
\begin{equation}\label{E135}
	\limsup_{y\to\infty}\frac{1}{\omega(y)}\int_0^y b(x,y)\omega(x)\,\rd x < 1.
\end{equation}
\hfill $\lozenge$
\end{remark}

\begin{remark}
The case when $\omega(x)=x^p$ with $p>1$ and $b$ satisfies the mass conservation condition \eqref{E41} 
is considered in \cite[Theorem~5.1.47]{banasiak_lamb_laurencot2020vol1}.
\hfill $\lozenge$
\end{remark}

In the following example we consider homogeneous fragmentation kernels;
these are also treated in \cite[Example~5.1.51]{banasiak_lamb_laurencot2020vol1}.

\begin{example}\label{E320}
Let
\[
	b(x,y) = \frac{1}{y}h\Bigl(\frac{x}{y}\Bigr)
\]
with a non-negative, measurable function $h$ such that $\int_0^1 h(\xi)\xi\,\rd\xi=1$.  
For instance, we can choose $h(z)=(\nu+2)z^\nu$ with $\nu\in(-2,0)$.
Then, for $\omega(x)=x^p$ with $p\ge1$, we have
\[
	\frac{1}{\omega(y)}\int_0^y b(x,y)\omega(x)\,\rd x
	= \int_0^y \frac{1}{y}h\Bigl(\frac{x}{y}\Bigr)\cdot\Bigl(\frac{x}{y}\Bigr)^p\,\rd x
	= \int_0^1 h(\xi)\xi^p\,\rd\xi
	\begin{cases} 
		= 1 & \text{if} \ p=1,
		\\[0.5ex]
		<1 & \text{if} \ p>1.
	\end{cases}
\]
This shows that \eqref{E41} is satisfied, i.e.\ we have mass conservation.
Moreover, for every $p>1$, \eqref{E213} holds with some $\kappa_2<1$ and any $\eta_0>0$,
and hence also \eqref{E135} is satisfied for every $p>1$.
\hfill $\lozenge$
\end{example}

\medskip

In the following proposition we compare the validity of \eqref{E135} for two different weights.

\begin{proposition}\label{E217}
Let $\omega_1$ and $\omega_2$ be differentiable and increasing functions such that $\omega_i(x)>0$ for $x>0$, $i\in\{1,2\}$.
Moreover, assume that
\begin{equation}\label{E218}
	\frac{\omega_1'(x)}{\omega_1(x)} \le \frac{\omega_2'(x)}{\omega_2(x)} \qquad \text{for} \ x>0.
\end{equation}
Then
\begin{equation}\label{E219}
	\limsup_{y\to\infty}\frac{1}{\omega_1(y)}\int_0^y b(x,y)\omega_1(x)\,\rd x
	\ge \limsup_{y\to\infty}\frac{1}{\omega_2(y)}\int_0^y b(x,y)\omega_2(x)\,\rd x.
\end{equation}
In particular, if \eqref{E135} holds with $\omega=\omega_1$, then it holds also with $\omega=\omega_2$.
\end{proposition}

\begin{proof}
Set $g_i(x)\coloneq\log\omega_i(x)$, $x>0$, $i\in\{1,2\}$.
For $0<x<y$ and $i\in\{1,2\}$ we can then write
\[
	\frac{\omega_i(x)}{\omega_i(y)} = e^{g_i(x)-g_i(y)}
	= \exp\biggl(-\int_x^y g_i'(\xi)\,\rd\xi\biggr)
	= \exp\biggl(-\int_x^y \frac{\omega_i'(\xi)}{\omega_i(\xi)}\,\rd\xi\biggr),
\]
which, together with \eqref{E218}, yields
\[
	\frac{1}{\omega_1(y)}\int_0^y b(x,y)\omega_1(x)\,\rd x
	\ge \frac{1}{\omega_2(y)}\int_0^y b(x,y)\omega_2(x)\,\rd x.
\]
This, in turn, implies \eqref{E219}.
\end{proof}

\section{Existence of a weight}\label{E13}

In this section we prove (in \cref{E127}) that, under the following mild assumption on the
fragmentation kernel $b$, there exists a weight $\omega$ such that \eqref{E212} and \eqref{E213} are satisfied.

\begin{assumption}\label{E143}
Assume that one of the following two statements is true:
\begin{myenum}
\item
	$b$ is bounded on $[0,\eta]^2$ for every $\eta>0$;
\item
	there exists $\eta_0>0$ such that $b$ is bounded on $[\eta_0,\eta]^2$ for every $\eta>\eta_0$
	and there exists $\omega_0\in L^\infty(0,\eta_0)$ such that $\omega_0(x)\ge0$ for almost every $x\in[0,\eta_0]$
	and $y\mapsto\int_0^{\eta_0} b(x,y)\omega_0(x)\,\rd x$ is locally bounded on $[\eta_0,\infty)$.
\end{myenum}
In case (i) we set $\eta_0=0$.
\hfill $\lozenge$
\end{assumption}

\begin{remark}\label{E142}
The assumption about the local boundedness of $y\mapsto\int_0^{\eta_0} b(x,y)\omega_0(x)\,\rd x$
is satisfied in many cases.  For instance, if $\int_0^y b(x,y)x\,\rd x\le y$ for all $y$, 
we can choose $\omega_0(x)=x$, which, for $y\in[\eta_0,\infty)$, yields
\[
	\int_0^{\eta_0} b(x,y)\omega_0(x)\,\rd x \le \int_0^y b(x,y)x\,\rd x \le y,
\]
which is locally bounded.
\hfill $\lozenge$
\end{remark}

\begin{theorem}\label{E127}
Suppose that \cref{E316,E143} hold and let $\kappa>0$.
Then there exists a function $\omega:[0,\infty)\to[0,\infty)$ such that $\omega(x)=\omega_0(x)$ for $x\in[0,\eta_0)$,
$\omega$ is continuous on $[\eta_0,\infty)$, $\omega(x)>0$ for $x>\eta_0$, and 
\[
	\int_0^y b(x,y)\omega(x)\,\rd x \le \kappa\omega(y) \qquad \text{for} \ y\in[\eta_0,\infty).
\]
\end{theorem}

\medskip

Before we prove \cref{E127}, we need two lemmas.

\begin{lemma}\label{E140}
Let $\eta_0$, $b$ and $\omega_0$ be as in \cref{E127}.
Then there exist continuous functions $h:[\eta_0,\infty)\to[0,\infty)$ and $\tilde b:[\eta_0,\infty)^2\to[0,\infty)$
such that
\begin{alignat}{2}
	\int_0^{\eta_0} b(x,y)\omega_0(x)\,\rd x &\le h(y), \qquad && y\in[\eta_0,\infty),
	\label{E138}
	\\[1ex]
	b(x,y) &\le \tilde b(x,y), \qquad && x,y\in[\eta_0,\infty).
	\label{E139}
\end{alignat}
\end{lemma}

\begin{proof}
Let us start with the proof of \eqref{E138}.  
By assumption, if we set
\[
	g(y) \coloneq \int_0^{\eta_0} b(x,y)\omega_0(x)\,\rd x, \qquad y\in[\eta_0,\infty),
\]
then it follows that $g$ is bounded on each interval $[\eta_0,\eta]$ with $\eta>\eta_0$,
and hence the numbers
\[
	h_n \coloneq \sup\bigl\{g(y):y\in[\eta_0,\eta_0+n+1]\bigr\}, \qquad n\in\NN_0,
\]
are well defined.  We construct $h$ as a piecewise linear function:
\[
	h(y) \coloneq h_n+(h_{n+1}-h_n)(y-\eta_0-n), \qquad \text{when} \ y\in[\eta_0+n,\eta_0+n+1).
\]
It is easy to see that $h$ is continuous, and that, for $y\in[\eta_0+n,\eta_0+n+1)$,
we have $h(y) \ge h_n \ge g(y)$, which proves \eqref{E138}.

The proof of \eqref{E139} is similar, but we construct a piecewise linear function $\tilde b$ which
is constant along diagonals.
Set
\[
	b_n \coloneq \sup\bigl\{b(x,y):x,y\in[\eta_0,\infty),\,x+y-2\eta_0\le n+1\bigr\}, \qquad n\in\NN_0,
\]
which is well defined since $b$ is bounded on $[\eta_0,\eta_0+n+1]^2$ by assumption.
It is easy to see that the piecewise linear function $\tilde b$, defined by
\begin{align*}
	& \tilde b(x,y) \coloneq b_n + (b_{n+1}-b_n)(x+y-2\eta_0-n) 
	\\[0.5ex]
	&\hspace*{30ex} \text{when} \ x,y\in[\eta_0,\infty),\,x+y-2\eta_0\in[n,n+1],
\end{align*}
is continuous on $[\eta_0,\infty)^2$.
Moreover, if $x,y\in[\eta_0,\infty)$ with $x+y-2\eta_0\in[n,n+1]$, 
then $\tilde b(x,y) \ge b_n \ge b(x,y)$.
\end{proof}

\medskip

The following lemma is based on standard ideas.  For the convenience of the reader we present the proof.

\begin{lemma}\label{E121}
Let $\eta_0\ge0$, let $\tilde b:[\eta_0,\infty)^2\to[0,\infty)$ and $f:[\eta_0,\infty)\to[0,\infty)$ be continuous functions,
and let $\kappa>0$.
Then there exists a unique continuous function $\omega:[\eta_0,\infty)\to[0,\infty)$ such that
\begin{equation}\label{E122}
	\int_{\eta_0}^y \tilde b(x,y)\omega(x)\,\rd x = \kappa\omega(y) - f(y), \qquad y\in[\eta_0,\infty).
\end{equation}
\end{lemma}

\begin{proof}
Let $\ell>0$ be arbitrary.  We construct a unique function $\omega_{[\ell]}$ on $[\eta_0,\eta_0+\ell]$ such that
\begin{equation}\label{E126}
	\int_{\eta_0}^y \tilde b(x,y)\omega_{[\ell]}(x)\,\rd x = \kappa\omega_{[\ell]}(y) - f(y), \qquad y\in[\eta_0,\eta_0+\ell].
\end{equation}
In the space $C([\eta_0,\eta_0+\ell])$ consider the positive Volterra operator $T_\ell$ defined by
\[
	(T_\ell u)(y) \coloneq \int_{\eta_0}^y \tilde b(x,y)u(x)\,\rd x, \qquad
	y\in[\eta_0,\eta_0+\ell],\; u\in C([\eta_0,\eta_0+\ell]).
\]
Set
\[
	M_\ell \coloneq \max_{x,y\in[\eta_0,\eta_0+\ell]}\tilde b(x,y).
\]
By induction we show that, for $u\in C([\eta_0,\eta_0+\ell])$,
\begin{equation}\label{E141}
	|(T_\ell^n u)(y)| \le \frac{1}{n!}M_\ell^n(y-\eta_0)^n\|u\|, \qquad
	y\in[\eta_0,\eta_0+\ell],\; n\in\NN_0,
\end{equation}
where $\|u\|\coloneq\max_{x\in[\eta_0,\eta_0+\ell]}|u(x)|$.
For $n=0$ this is trivial.  Now assume that \eqref{E141} holds for some $n$; then, for $y\in[\eta_0,\eta_0+\ell]$,
\begin{align*}
	|(T_\ell^{n+1}u)(y)| &\le \int_{\eta_0}^y \tilde b(x,y)|(T_\ell^n u)(x)|\,\rd x
	\le \int_{\eta_0}^y M_\ell\frac{1}{n!}M_\ell^n(x-\eta_0)^n\|u\|\,\rd x
	\\[1ex]
	&= \frac{1}{(n+1)!}M_\ell^{n+1}(y-\eta_0)^{n+1}\|u\|.
\end{align*}
Hence \eqref{E141} is true for all $n\in\NN_0$, which implies that
\[
	\|T_\ell^n\| \le \frac{M_\ell^n\ell^n}{n!}, \qquad n\in\NN_0.
\]
It follows that
\[
	R(\kappa,T_\ell) = \frac{1}{\kappa}\sum_{n=0}^\infty \frac{1}{\kappa^n}T_\ell^n,
\]
where the series converges in the operator norm, and therefore $R(\kappa,T_\ell)$ is a positive operator.
This shows that \eqref{E126} has a unique positive solution.
Since $\ell$ was arbitrary and $\omega_{[\ell_2]}$ is an extension of $\omega_{[\ell_1]}$ if $\ell_1<\ell_2$,
we obtain a unique positive solution of \eqref{E122}.
\end{proof}

\begin{proof}[Proof of \cref{E127}]
By \cref{E140} there exist $h$ and $\tilde b$ such that \eqref{E138} and \eqref{E139} hold.
Without loss of generality, $h$ can be chosen such that $h(y)>0$ for $y>\eta_0$. 
Let $\omega$ be the unique solution of \eqref{E122} with $f \coloneq h$ and set $\omega(x)\coloneq \omega_0(x)$
for $x\in[0,\eta_0)$.  Then, for $y\in[\eta_0,\infty)$,
\begin{align*}
	\int_0^y b(x,y)\omega(x)\,\rd x
	&= \int_0^{\eta_0} b(x,y)\omega_0(x)\,\rd x + \int_{\eta_0}^y b(x,y)\omega(x)\,\rd x
	\\[1ex]
	&\le h(y) + \int_{\eta_0}^y \tilde b(x,y)\omega(x)\,\rd x
	= \kappa\omega(y).
\end{align*}
Since $h(y)>0$ for $y>\eta_0$, we have $\omega(y)>0$ for $y>\eta_0$.
\end{proof}

\medskip

The next theorem shows that, under certain assumptions, condition \eqref{E135} is
satisfied with an exponential weight.

\begin{theorem}\label{E132}
Assume that
\begin{equation}\label{E133}
  \int_0^y b(x,y)x\,\rd x \le y \qquad \text{for all} \ y>0
\end{equation}
and that there exist $\delta_1,\delta_2>0$, $d>1$ and $\bm>0$ such that
\begin{align}
  \int_0^{\delta_1} b(x,y)\,\rd x &\le d^y,
  \label{E137}
  \\[1ex]
  b(x,y) &\le \bm, \qquad x\in[y-\delta_2,y],
  \label{E134}
\end{align}
for large enough $y$.  Then there exists $c>1$ such that \eqref{E135}
is satisfied with $\omega(x)=c^x$.
\end{theorem}

\begin{proof}
Let $c>d$ be such that $\frac{1}{\log c}<\delta_1$ but otherwise arbitrary, and let $\delta\in(0,\delta_2]$.
Note that the function $x\mapsto\frac{c^x}{x}$ is strictly increasing 
on the interval $\bigl[\frac{1}{\log c},\infty\bigr)$, which, together 
with \eqref{E133}, \eqref{E134} and \eqref{E137}, implies that
\begin{align*}
	\int_0^y b(x,y)\omega(x)\,\rd x 
	&= \int_0^{\frac{1}{\log c}} b(x,y)c^x\,\rd x + \int_{\frac{1}{\log c}}^{y-\delta} b(x,y)\frac{c^x}{x}x\,\rd x
	+ \int_{y-\delta}^y b(x,y)c^x\,\rd x
	\\[1ex]
	&\le c^{\frac{1}{\log c}}\int_0^{\frac{1}{\log c}} b(x,y)\,\rd x 
	+ \frac{c^{y-\delta}}{y-\delta}\int_{\frac{1}{\log c}}^{y-\delta} b(x,y)x\,\rd x
	+ \delta \bm c^y
	\\[1ex]
	&\le e\int_0^{\delta_1} b(x,y)\,\rd x
	+ \frac{c^{y-\delta}}{y-\delta}\int_0^y b(x,y)x\,\rd x + \delta \bm c^y
	\\[1ex]
	&\le ed^y + \frac{y}{y-\delta}c^{y-\delta} + \delta \bm c^y.
\end{align*}
Since $c>d$, we arrive at
\begin{equation}\label{E136}
  \limsup_{y\to\infty}\frac{1}{c^y}\int_0^y b(x,y)c^x\,\rd x
  \le \limsup_{y\to\infty}\biggl[e\Bigl(\frac{d}{c}\Bigr)^y + \frac{y}{y-\delta}c^{-\delta} + \delta \bm\biggr]
  = c^{-\delta} + \delta \bm.
\end{equation}
If we first choose $\delta$ small enough so that $\delta \bm<\frac12$ and then $c$ large enough so that $c^{-\delta}<\frac12$, 
then the right-hand side of \eqref{E136} is strictly less than 1.
\end{proof}

The theorem can be applied to show the existence of an exponential weight in \cref{E130}.

\section{Examples}\label{E14}

Let us consider two examples.  The first one shows that \eqref{E135} can be satisfied for an exponential weight
although it is not satisfied for any power weight.  It yields binary fragmentation and 
is similar to the Becker--D\"oring model in the discrete case 
in the sense that particles of size $y$ fragment into small particles (of size at most $1$) 
and large particles (of size at least $y-1$).

\begin{example}\label{E130}
Let
\[
	b(x,y) 
	= \begin{cases}
		1, & x\in[0,1]\cup[y-1,y], \\[0.5ex]
		0, & x\in(1,y-1),
	\end{cases}
\]
when $y>2$ and $b(x,y)=\frac{2}{y}$ when $y\le2$.
Then we have mass conservation: e.g.\ for $y>2$ we have
\[
	\int_0^y b(x,y)x\,\rd x = \int_0^1 x\,\rd x + \int_{y-1}^y x\,\rd x
	= \frac{1}{2}\Bigl(1+y^2-(y-1)^2\Bigr) = y.
\]
For the power weight $\omega(x)=x^p$ with $p\ge1$ we obtain from the Mean Value Theorem that,
for $y>2$,
\begin{align*}
	\int_0^y b(x,y)\omega(x)\,\rd x &= \int_0^1 x^p\,\rd x + \int_{y-1}^y x^p\,\rd x
	= \frac{1}{p+1}\Bigl(1+y^{p+1}-(y-1)^{p+1}\Bigr)
	\\[1ex]
	&= \frac{1}{p+1}\Bigl(1+(p+1)\xi_y^p\Bigr) \hspace*{10ex} \text{(with $\xi_y\in(y-1,y)$)}
	\\[1ex]
	&\sim y^p = \omega(y)
\end{align*}
as $y\to\infty$.  Hence \eqref{E135} is not satisfied.
On the other hand, we obtain from \cref{E132} that there exists an exponential weight so that \eqref{E135} is fulfilled.
We can also give a direct and explicit construction, namely, for the exponential weight $\omega(x)=e^x$, we have,
for $y>2$,
\begin{align*}
	\frac{1}{\omega(y)}\int_0^y b(x,y)\omega(x)\,\rd x 
	&= \frac{1}{e^y}\biggl[\int_0^1 e^x\,\rd x + \int_{y-1}^y e^x\,\rd x\biggr]
	\\[1ex]
	&= \frac{1}{e^y}\Bigl(e-1+e^y-e^{y-1}\Bigr)
	\to 1-\frac{1}{e}
\end{align*}
as $y\to\infty$, and hence \eqref{E135} is satisfied.
\hfill $\lozenge$
\end{example}

The second example, which is from \cite[Example~5.1.51]{banasiak_lamb_laurencot2020vol1},
shows that \eqref{E135} is not always satisfied for an exponential weight.
This is different from the situation in the discrete case;
for the latter see \cite[Theorem~5.5]{kerr_lamb_langer2020fragpaper}.

\begin{example}\label{E131}
Let
\[
	b(x,y) 
	= \begin{cases}
		y, & x\in\bigl[0,\frac{1}{y}\bigr]\cup\bigl[y-\frac{1}{y},y\bigr], \\[1.5ex]
		0, & x\in\bigl(\frac{1}{y},y-\frac{1}{y}\bigr),
	\end{cases}
\]
when $y>\sqrt{2}$ and $b(x,y)=\frac{2}{y}$ when $y\le\sqrt{2}$.
For $\omega(x)=e^x$ we have, for $y>\sqrt{2}$, (big O notation for $y\to\infty$)
\begin{align*}
	\int_0^y b(x,y)\omega(x)\,\rd x 
	&= y\int_0^{\frac{1}{y}}e^x\,\rd x + y\int_{y-\frac{1}{y}}^y e^x\,\rd x
	\\[1ex]
	&= y\bigl(e^{\frac{1}{y}}-1\bigr) + y\bigl(e^y-e^{y-\frac{1}{y}}\bigr)
	= y\bigl(e^{\frac{1}{y}}-1\bigr) + ye^y\bigl(1-e^{-\frac{1}{y}}\bigr)
	\\[1ex]
	&= y\biggl[1+\frac{1}{y}+\BigO\Bigl(\frac{1}{y^2}\Bigr)-1\biggr]
	+ ye^y\biggl[1-\biggl(1-\frac{1}{y}+\BigO\Bigl(\frac{1}{y^2}\Bigr)\biggr)\biggr]
	\\[1ex]
	&= 1+\BigO\Bigl(\frac{1}{y}\Bigr) + e^y\biggl(1+\BigO\Bigl(\frac{1}{y}\Bigr)\biggr)
	\sim e^y.
\end{align*}
Hence \eqref{E135} is not satisfied.  A similar---but slightly lengthier---calculation can 
be done for an arbitrary exponential weight $\omega(x)=c^x$ with $c>1$.
As shown in \cite[Example~5.1.51]{banasiak_lamb_laurencot2020vol1}, condition \eqref{E135} 
is not satisfied for power weights either.
However, we can find a faster growing weight so that \eqref{E135} holds, 
namely, choose $\omega(x)=xe^{x^2}$.  Then, for $y>\sqrt{2}$,
\begin{align*}
	\frac{1}{\omega(y)}\int_0^y b(x,y)\omega(x)\,\rd x
	&= \frac{1}{ye^{y^2}}\biggl[\int_0^{\frac{1}{y}}yxe^{x^2}\,\rd x+\int_{y-\frac{1}{y}}^y yxe^{x^2}\,\rd x\biggr]
	\\[1ex]
	&= e^{-y^2}\frac{1}{2}\biggl[e^{\frac{1}{y^2}}-1+e^{y^2}-e^{(y-\frac{1}{y})^2}\biggr]
	\\[1ex]
	&\sim \frac{1}{2}\Bigl(1-e^{-2+\frac{1}{y^2}}\Bigr)
	\to \frac{1-e^{-2}}{2} < 1
\end{align*}
as $y\to\infty$, which shows that \eqref{E135} is satisfied.
\hfill $\lozenge$
\end{example}


\noindent
Lyndsay Kerr, Wilson Lamb, Matthias Langer \\[1ex]
Department of Mathematics and Statistics \\
University of Strathclyde \\
26 Richmond Street \\
Glasgow G1 1XH \\
United Kingdom \\[1ex]
Email addresses: \\
\texttt{lyndsay.kerr@strath.ac.uk} \\ 
\texttt{w.lamb@strath.ac.uk} \\
\texttt{m.langer@strath.ac.uk} \\[1ex]
\texttt{ORCID:} \\
0000-0002-6667-7175 (L.K.) \\
0000-0001-8084-6054 (W.L.) \\
0000-0001-8813-7914 (M.L.)

\end{document}